\documentclass[a4paper,11pt]{article}%
\usepackage{graphicx}
\usepackage{geometry}
\usepackage{amsmath}
\usepackage{amssymb}
\usepackage{amsfonts}
\usepackage{amsthm,amscd}%
\providecommand{\U}[1]{\protect\rule{.1in}{.1in}}
\def\figurename{Figure}
\makeatletter
\renewcommand{\fnum@figure}[1]{\figurename~\thefigure.}
\makeatother
\def\tablename{Table}
\makeatletter
\renewcommand{\fnum@table}[1]{\tablename~\thetable.}
\makeatother
\def \bop {\noindent\textbf{Proof. }}
\def \eop {\hbox{}\nobreak\hfill
\vrule width 2mm height 2mm depth 0mm
\par \goodbreak \smallskip}
\newtheorem{theorem}{Theorem}[section]
\newtheorem{lemma}[theorem]{Lemma}
\newtheorem{corollary}[theorem]{Corollary}
\newtheorem{proposition}[theorem]{Proposition}
\theoremstyle{definition}
\newtheorem{definition}[theorem]{Definition}

\theoremstyle{remark}
\newtheorem{remark}[theorem]{Remark}
\numberwithin{equation}{section}

\begin{document}

\title{Approximation and generic properties of McKean-Vlasov stochastic equations
with continuous coefficients}
\author{{Mohamed Amine Mezerdi }\thanks{ Laboratory of Applied Mathematics, University
of Biskra, Po. Box 145, Biskra (07000), Algeria. \textit{(E-mail:
amine.mezerdi@univ-biskra.dz)}}
\and {Khaled Bahlali}\thanks{ Laboratoire IMATH, Universit\'{e} du Sud-Toulon-Var,
B.P 20132, 83957 La Garde Cedex 05, France. \textit{(E-mail:
bahlali@univ-tln.fr)}}
\and {Nabil Khelfallah }\thanks{ Laboratory of Applied Mathematics, University of
Biskra, Po. Box 145, Biskra (07000), Algeria. \textit{(E-mail:
nabilkhelfallah@yahoo.fr)}}
\and Brahim Mezerdi\thanks{ King Fahd University of Petroleum and Minerals,
Department of Mathematics and Statistics, P.O. Box 1916, Dhahran 31261, Saudi
Arabia. (\textit{E-mail: brahim.mezerdi@kfupm.edu.sa)}}}
\maketitle
\date{}

\begin{abstract}
We consider various approximation properties for systems driven by a Mc
Kean-Vlasov stochastic differential equations (MVSDEs) with continuous
coefficients, for which pathwise uniqueness holds. We prove that the solution
of such equations is stable with respect to small perturbation of initial
conditions, parameters and driving processes. Moreover, the unique strong
solutions may be constructed by effective approximation procedures, without
using the famous Yamada-Watanabe theorem. Finally we show that the set of
bounded uniformly continuous coefficients for which the corresponding MVSDE
have a unique strong solution is a set of second category in the sense of Baire.



\end{abstract}

\section{Introduction}

McKean-Vlasov stochastic differential equations (MVSDE) have been investigated
by McKean \cite{McK}, for the first time, as the counterpart of Vlasov
\cite{Vla} non linear partial differential equations (PDE) arising in
statistical physics. They describe the limiting bahaviour of an individual
particle evolving within a large system of particles, with weak interaction,
as the number of particles tends to infinity. These equations are called non
linear SDEs in the sense that the coefficients depend not only on the state
variable, but also on its marginal distribution and their solutions are called
non linear diffusions. A pedagogical and rigorous treatment of these equations
appear in the seminal Saint-Flour course by Sznitman \cite{Sn}$.$ In the last
decade, there has been a renewed interest for this kind of equations, due to
their intimate relationship with the so-called mean-field games, introduced
independently by Lasry-Lions \cite{LasLio} and Huang, Malham\'{e} and
Caines$\ $\cite{HMC} to find approximate Nash equilibriums for differential
games with a large number of players$.$ For a complete and detailed treatment
of mean-field games and mean-field control problems, one can refer to the
excellent book by Carmona and Delarue \cite{CarDel} and its complete list of
references$.$ Backward stochastic differential equations of the McKean-Vlasov
type have been investigated in particular in \cite{BDLP, BLP} just to cite a
few references.

Since the pioneering work of McKean \cite{McK}, a huge literature on
existence, uniqueness, numerical schemes and propagation of chaos theorems was
developped. Existence and uniqueness of strong solutions were obtained under
global Lipschitz coefficients in \cite{Gra, JMW, Sn} by using the fixed point
theorem on the space of continuous functions with values on the space of
probability measures, equipped with Wasserstein distance. MVSDEs with non
regular coefficients appear naturally in many mean-field models. The so-called
mean-field FitzHugh-Nagumo model and the network of Hodgkin-Huxley neurons are
typical examples (see \cite{BFT, DEG}). It is clear that if the coefficients
are not globally Lipschitz, the Gronwall inequality and its variants fail, so
that fixed point theorems are no longer applicable. It is well known that for
It\^{o} SDEs with singular drift, the Brownian motion brings a regularization
effect, provided the diffusion matrix is uniformly nondegenerate (see
\cite{ZK} for It\^{o} SDEs). Recently, many existence, uniqueness results were
proved for non degenerate MSVSDEs under non Lipschitz coefficients \cite{Ch,
ChFr, Chia, HSS, MiVe}. For some recent results on the numerical aspects see
\cite{DEG, GP}. Let us point out that contrary to It\^{o}'s SDEs, regularity
assumptions of local nature on the coefficients, such as locally Lipschitz
coefficients do not lead to unique (local) strong solutions (see \cite{Sheu}
for counterexamples).

It is a well known that if the coefficients are Lipschitz continuous, then
MVSDE (\ref{MVSDE}) has a unique strong solution $X_{t}\left(  x\right)  $,
which is continuous with respect to the initial condition and coefficients.
Moreover, the solution may be constructed by means of various numerical
schemes (see \cite{BMMA}).

Our purpose in this paper, is to study strong stability properties of the
solution of (\ref{MVSDE}) under pathwise uniqueness of solutions and merely
continuous coefficients. Since the coefficients are only continuous without
additional regularity, one cannot expect to apply Gronwall's lemma. Instead of
Gronwall's lemma, we use tightness arguments and the famous Skorokhod
selection theorem to prove the desired convergence results. Of course we
should not expect precise convergence speed as this last property is based on
regularity of the coefficients.

The paper is organized as follows. In the second section we prove that the
Euler polygonal scheme is convergent provided that there is pathwise
uniqueness. This provides us with an effective way to construct strong
solutions for MVSDEs without appealing to the famous Yamada-Watanabe
theorem.\ In particular this results extends \cite{BMO1, GK, KN} proved for
It\^{o} SDEs to MVSDEs. In the third section we prove that the solution is
stable under small perturbion of to initial condition, coefficients and
driving processes. In the last section, we show that the set of bounded
uniformly continuous coefficients for which strong existence and uniqueness
hold is a generic property in the sense of Baire. This means that in the sense
of Baire category, most of MVSDEs with bounded uniformly continuous
coefficients have unique solutions. \ This last result extends in particular
\cite{BMO1, BMO2} to MVSDEs.

\section{Assumptions and preliminaries}

\subsection{The Wasserstein distance}

\begin{definition}
Let $(M,d)$ be a metric space, for which every probability measure on $M$ is a
Radon measure (a so-called Radon space). Denote $\mathcal{P}_{p}(M)$ the
collection of all probability measures $\mu$ on $M$ with finite moment of
order $p$ for some $x_{0}$ in $M$, $%
{\displaystyle\int\nolimits_{M}}
d(x_{0},x)^{p}\mu(dx)<+\infty.$ Then the $p-$Wasserstein distance between two
probability measures $\mu$ and $\nu$ in $\mathcal{P}_{p}(M)$ is defined as
\end{definition}

\begin{center}
$W_{p}$($\mu$,$\nu)=\left(  \inf\left\{
{\displaystyle\int\nolimits_{M\times M}}
d(x,y)^{p}\gamma(dx,dy);\text{ }\gamma\in\Gamma(\mu,\nu)\right\}  \right)
^{1/p}$
\end{center}

where $\Gamma(\mu,\nu)$ denotes the collection of all measures on $M\times M$
with marginals $\mu$ and $\nu$ on the first and second factors respectively.
The set $\Gamma(\mu,\nu)$ is also called the set of all couplings of $\mu$ and
$\nu$

\smallskip

The Wasserstein metric may be equivalently defined by $W_{p}$($\mu$%
,$\nu)=\left(  \inf E\left[  d(X,Y)^{p}\right]  \right)  ^{1/p}$ where the
infimum is taken over all the joint probability distributions of the random
variables $X$ and $Y$ with marginals $\mu$ and $\nu.$

In the case where the metric space is replaced by the euclidian space
$\mathbb{R}^{d},$ then the $p$-Wasserstein distance $W_{p}(\mu,\nu)$ is
defined by:
\[
W_{p}(\mu,\nu)^{p}=\inf_{\gamma\in\Gamma(\mu,\nu)}\int_{\mathbb{R}^{d}%
\times\mathbb{R}^{d}}\left\vert x-y\right\vert ^{p}d\pi(x,y)^{1/p}%
\]

or equivalently $W_{p}(\mu,\nu)^{p}=\inf\mathbb{E[}\left\vert X-Y\right\vert
^{p}].$

In particular if $X$ and $Y$\ are square integrable random variables, we have
$W_{2}(P_{X},P_{Y})\leq\mathbb{E[}\left\vert X-Y\right\vert ^{2}]^{1/2}.$

In the literature the Wasserstein metric is restricted to $W_{2}$ while
$W_{1}$ is often called the Kantorovich-Rubinstein distance because of the
role it plays in optimal transport.

\subsection{Assumptions}

Let $\left(  B_{t}\right)  $ a $d$-dimensional Brownian motion defined on a
probablity space $(\Omega,\mathcal{F},P),$ equipped with a filtration $\left(
\mathcal{F}_{t}\right)  ,$ satisfying the usual conditions. Throughout this
paper, we consider McKean-Vlasov stochastic differential equation (MVSDE),
called also mean-field stochastic differential equation of the form

\begin{center}%
\begin{equation}
\left\{
\begin{array}
[c]{l}%
dX_{t}=b(t,X_{t},\mathbb{P}_{X_{t}})dt+\sigma(t,X_{t},\mathbb{P}_{X_{t}%
})dB_{t}\\
X_{0}=x
\end{array}
\right.  \label{MVSDE}%
\end{equation}

\end{center}

For this kind of stochastic differential equations, the drift $b$ and
diffusion coefficient $\sigma$ depend not only on the state process $X_{t}$,
but also on its marginal distribution $\mathbb{P}_{X_{t}}$.

Assume that the coefficients satisfy the following conditions.

(\textbf{H}$_{\mathbf{1}}$\textbf{)} Assume that

\begin{center}
$%
\begin{array}
[c]{c}%
b:[0,T]\times\mathbb{R}^{d}\times\mathcal{P}_{2}(\mathbb{R}^{d}%
)\longrightarrow\mathbb{R}^{d}\\
\sigma:[0,T]\times\mathbb{R}^{d}\times\mathcal{P}_{2}(\mathbb{R}%
^{d})\longrightarrow\mathbb{R}^{d}\otimes\mathbb{R}^{d}%
\end{array}
$
\end{center}

are Borel measurable functions and continuous in $(x,\mu)$ uniformly in
$t\in\left[  0,T\right]  .$

(\textbf{H}$_{\mathbf{2}}$\textbf{)} There exist $C>0$ such that for any
$t\in\lbrack0,T],x\in\mathbb{R}^{d}$ and $\mu\in\mathcal{P}_{2}(\mathbb{R}%
^{d}),$

\begin{center}
$|b(t,x,\mu)|\leq C\left(  1+\left\vert x\right\vert +W_{2}(\mu,\delta
_{0}\right)  ),$

$|\sigma(t,x,\mu)|\leq C\left(  1+\left\vert x\right\vert +W_{2}(\mu
,\delta_{0}\right)  ),$
\end{center}

where $W_{2}$ is the 2-Wasserstein distance and $\delta_{0}$ is the Dirac
measure at $0$.

The following theorem states that under global Lipschitz condition,
(\ref{MVSDE}) admits a unique solution. Its complete proof is given in
\cite{Sn} for a drift depending linearly on the law of $X_{t}$ that is
$b(t,x,\mu)=%
{\displaystyle\int\limits_{\mathbb{R}^{d}}}
b^{\prime}(t,x,y)\mu(dy)$ and a constant diffusion. The general case as
(\ref{MVSDE}) is treated in \cite{CarDel} Theorem 4.21 or \cite{JMW}
Proposition 1.2.\ The proof is based on a fixed point theorem on the space of
continuous functions with values in $\mathcal{P}_{2}(\mathbb{R}^{d})$ endowed
with Wasserstein metric. Note that in \cite{Gra, JMW} the authors consider
MVSDEs driven by general L\'{e}vy process instead of a Brownian motion.

\begin{theorem}
Assume $\mathbf{(H_{1})}$, $\mathbf{(H_{2})}$ and

$(\mathbf{H}_{\mathbf{3}})$ there exist $L>0$ such that for any $t\in
\lbrack0,T],x,$ $x^{\prime}\in\mathbb{R}^{d}$ and $\mu,$ $\mu^{\prime}%
\in\mathcal{P}_{2}(\mathbb{R}^{d}),$

$|b(t,x,\mu)-b(t,x^{\prime},\mu^{\prime})|\leq C\left(  \left\vert
x-x^{\prime}\right\vert +W_{2}(\mu,\mu^{\prime}\right)  ),$

$|\sigma(t,x,\mu)-\sigma(t,x^{\prime},\mu^{\prime})|\leq C\left(  \left\vert
x-x^{\prime}\right\vert +W_{2}(\mu,\mu^{\prime}\right)  ),$

then MVSDE (\ref{MVSDE}) admits a unique solution such that $E[\sup_{t\leq
T}|X_{t}|^{2}]<+\infty.$
\end{theorem}

\bop See \cite{CarDel} Theorem 4.21$.$\eop

Other versions of the MVSDEs, which are particular cases of (\ref{MVSDE}) have
been considered in the literature.

1) The following MVSDE has been treated in literature%

\begin{equation}
\left\{
\begin{array}
[c]{l}%
dX_{t}=b(t,X_{t},%
{\displaystyle\int}
\varphi(y)\mathbb{P}_{X_{t}}(dy))dt+\sigma(t,X_{t},%
{\displaystyle\int}
\psi(y)\mathbb{P}_{X_{t}}(dy))dW_{t}\\
X_{0}=x,
\end{array}
\right.  \label{MVSDE1}%
\end{equation}

2) MVSDEs studied in the framework of statistical physics take the form

\begin{center}
$\left\{
\begin{array}
[c]{l}%
dX_{t}=%
{\displaystyle\int\limits_{\mathbb{R}^{d}}}
b(t,X_{t},y)\mathbb{P}_{X_{t}}(dy)dt+dB_{t}\\
X_{0}=x
\end{array}
\right.  $
\end{center}

where $b:[0,T]\times\mathbb{R}^{d}\times\mathbb{R}^{d}\longrightarrow
\mathbb{R}^{d}$ is a Borel measurable function such that $b(t,.,.)$ is
Lipschitz. This is a particular class of MVSDEs for interacting diffusions,
considered by McKean \ (see \cite{Sn} for details), where the drift is linear
on the probability distribution.\ It is easy to see that the drift is
Lipschitz in the measure variable with respect to Wasserstein metric.

The definition of pathwise uniqueness for equation (\ref{MVSDE}) is given by
the following.

\begin{definition}
\textit{\ }We say that pathwise uniqueness holds for equation (\ref{MVSDE}%
)\ if \ $X$ and $X^{\prime}$ are two solutions defined on the same probability
space $\left(  \Omega,\mathcal{F}\text{,}P\right)  $ with common Brownian
motion $\left(  B\right)  $, with possibly different filtrations such that
$P\left[  X_{0}=X_{0}^{\prime}\right]  =1$, then $X$\ and $X^{\prime}$\ are indistinguishable.
\end{definition}

Let us recall Kolmogorov's tightness criteia for stochastic processes and
Skorokhod selection theorem, which will be extensively used in the sequel.

\begin{lemma}
(Skorokhod selection theorem $\cite{IW}$ page 9) \textit{Let }$\left(
S,\rho\right)  $\textit{\ be a complete separable metric space, }%
$P_{n},n=1,2,...$ \textit{and }$P$\textit{\ be probability measures on
}$\left(  S,\mathcal{B}\left(  S\right)  \right)  $ \textit{such that }$P_{n}%
$\textit{\ }$\underset{n\rightarrow+\infty}{\longrightarrow}P$.\textit{\ Then,
on a probability space }$\left(  \widehat{\Omega},\widehat{\mathcal{F}%
},\widehat{P}\right)  $, \textit{we can construct S-valued random variables
}$X_{n}$\textit{, }$n=1,2,...,$ \textit{and }$X$\textit{\ such that:}
\end{lemma}

\textit{(i) }$P_{n}=\widehat{P}^{X_{n}},n=1,2,...$, \textit{and }%
$P=\widehat{P}^{X} $.

\textit{(ii) }$X_{n}$\textit{\ converges to }$X,$\textit{\ }$\widehat{P}%
$\textit{\ almost surely.}

\smallskip\ 

\begin{lemma}
(Skorokhod limit theorem $\cite{Sk}$) Let $\left(  \Omega,\mathcal{F}%
\text{,}P\right)  $ be a probability space, $(H_{n})$ be uniformly bounded
processes and $(H^{n})$ be a sequence of Brownian motions defined on the same
space such that the stochastic integral $%
{\displaystyle\int\limits_{O}^{T}}
H_{s}^{n}dW_{s}^{n}$ are well defined for each $n\geq0.$ Assume moreover that

a) $\lim_{h\rightarrow0}\sup_{n}\sup_{\left\vert s-t\right\vert <h}P\left(
\left\vert H_{s}^{n}-H_{t}^{n}\right\vert >\varepsilon\right)  =0$

b) $\left(  H_{s}^{n},W_{s}^{n}\right)  $ converges to $\left(  H_{s}%
^{0},W_{s}^{0}\right)  $ in probability.

Then $%
{\displaystyle\int\limits_{O}^{T}}
H_{s}^{n}dW_{s}^{n}$ converges in probability to $%
{\displaystyle\int\limits_{O}^{T}}
H_{s}^{0}dW_{s}^{0}$
\end{lemma}

\begin{lemma}
$\left(  \text{Kolmogorov criterion for tightness }\cite{IW}\text{ page
18}\right)  $ \textit{Let }$\left(  X_{n}\left(  t\right)  \right)  $,
$n=1,2,...$\textit{, be a sequence of d-dimensional continuous processes
satisfying the following two conditions:}
\end{lemma}

\textit{(i) There exist positive constants }$M$\textit{\ and }$\gamma
$\textit{\ such that }$E\left[  \left|  X_{n}\left(  0\right)  \right|
^{\gamma}\right]  \leq M$\textit{\ for every }$n=1,2,.....$

\textit{(ii) There exist positive constants }$\alpha,$ $\beta,$ $M_{k},$
$k=1,2,...,$ \textit{such that:\newline}$E\left[  \left|  X_{n}\left(
t\right)  -X_{n}\left(  s\right)  \right|  ^{\alpha}\right]  \leq M_{k}\left|
t-s\right|  ^{1+\beta}$\textit{\ for every }$n$\textit{\ and }$t,s\in\left[
0,k\right]  ,(k=1,2,...).$

\textit{Then there exist a subsequence }$(n_{k})$\textit{, a probability space
}$\left(  \widehat{\Omega},\widehat{\mathcal{F}},\widehat{P}\right)  $
\textit{and d-dimensional continuous processes }$\widehat{X}_{n_{k}%
},k=1,2,..., $ \textit{and }$\widehat{X}$\textit{\ defined on it such that }

\textit{1) The laws of }$\widehat{X}_{n_{k}}$\textit{\ and }$X_{n_{k}}%
$\textit{\ coincide.}

\textit{2) }$\widehat{X}_{n_{k}}(t)$\textit{\ converges to }$\widehat{X}%
(t)$\textit{\ uniformly on every finite time interval }$\widehat{P}%
$\textit{\ almost surely}.

\section{Construction of strong solutions by approximation}

It is well known for classical It\^{o} SDEs \cite{IW} as well as for
McKeanVlasov SDEs \cite{Kur} that weak existence and pathwise uniqueness imply
the existence and uniqueness of a strong solution.\ This is a corollary of the
famous Yamada-Watanabe theorem (see \cite{IW}). In this section we prove that
under the pathwise uniqueness, the strong solution, may be constructed by
means of an approximation procedure and may be written as a measurable
functional of the initial condition and the Brownian motion without appealing
to the famous Yamada-Watanabe theorem.

Let $\left(  \Delta^{n}\right)  $ be a sequence of partitions of the interval
$\left[  0,T\right]  $ where $\Delta^{n}:$ $0=t_{0}^{n}<t_{1}^{n}%
<.....<t_{n}^{n}=T$ such that

\begin{center}
$\lim_{n\rightarrow+\infty}\left\Vert \Delta^{n}\right\Vert =\lim
_{n\rightarrow+\infty}\underset{i}{\max}\left(  t_{i+1}^{n}-t_{i}^{n}\right)
=0$
\end{center}

Define the Euler polygonal approximation for equation (\ref{MVSDE}) by:

\begin{center}
$X_{\Delta^{n}}(x,t)=x+%
{\displaystyle\int_{0}^{t}}
b(\phi_{\Delta^{n}}(s),X_{\Delta^{n}},P_{\Delta^{n}})ds+%
{\displaystyle\int_{0}^{t}}
\sigma(\phi_{\Delta^{n}}(s),X_{\Delta^{n}},P_{X_{\Delta^{n}}})dB_{s}$
\end{center}

where $\phi_{\Delta^{n}}(s)=t_{i},$ if $t_{i}^{n}\leq s<t_{i+1}^{n}$ and
$\left\Vert \Delta^{n}\right\Vert =\underset{i}{\max}\left(  t_{i+1}^{n}%
-t_{i}^{n}\right)  $ and $X_{\Delta^{n}}=X_{\Delta^{n}}(x,\phi_{\Delta^{n}%
}(s))$

\begin{theorem}
Assume $\left(  \mathbf{H}_{\mathbf{1}}\right)  $ and $\left(  \mathbf{H}%
_{\mathbf{2}}\right)  $, then under pathwise uniqueness we have:

1)$\underset{n\rightarrow0}{\lim}\underset{}{E\left[  \underset{t\leq T}{\sup
}\left\vert X_{\Delta^{n}}(x,t)-X(x,t)\right\vert ^{2}\right]  }=0$

2) There exists a measurable functional $F:\mathbb{R}^{d}\times W_{0}%
^{d}\longrightarrow W^{d}$ which is adapted such that the unique solution
$X_{t}$ can be written $X(.)=F(X(0),B(.)),$ where $W^{d}=C\left(
\mathbb{R}_{+},\mathbb{R}^{d}\right)  $ and $W_{0}^{d}=\left\{  w\in C\left(
\mathbb{R}_{+},\mathbb{R}^{d}\right)  :w(0)=0\right\}  $ are equipped with
their Borel $\sigma-$fields and the filtrations of coordinates.
\end{theorem}

\bop1) Suppose that the conclusion of our theorem is false, then there exists
a sequence$\left(  \Delta_{n}\right)  $ and $\delta\geq0$ such that%

\begin{equation}
\underset{n\rightarrow\infty}{\lim\inf}E\left[  \underset{t\leq T}{\sup
}\left\vert X_{\Delta_{n}}(x,t)-X_{t}\right\vert ^{2}\right]  \geq\delta.
\label{Hyp1}%
\end{equation}

Let $\mathcal{C}$($\left[  0,T\right]  )$ be the space of continuous functions
equipped with the topology of uniform convergence and $\mathcal{P}_{2}\left(
\mathcal{C}(\left[  0,T\right]  )\right)  $ the space of probability measures
equipped with the Wasserstein metric.

Using assumptions $\mathbf{(H_{1})}$, $\mathbf{(H_{2})}$ and classical
arguments of stochastic calculus, it is easy to see that the sequence
$(X_{\Delta_{n}},X,B,P_{X_{\Delta_{n}}},P_{X})$ satisfies the conditions of
Kolmogorov criteria, then it is tight in $\mathcal{C}$($\left[  0,T\right]
)^{3}\times\mathcal{P}_{2}\left(  \mathcal{C}(\left[  0,T\right]  )\right)
^{2}.$

Then by Skorokhod limitTheorem (See Appendix), there exist a probability space
$\left(  \widehat{\Omega},\widehat{\mathcal{F}},\widehat{P}\right)  $ and a
sequence of stochastic processes $\left(  \widehat{X_{t}^{n}},\widehat
{Y_{t}^{n}},\widehat{B_{t}^{n}},\widehat{\mu_{t}^{n}},\widehat{\nu_{t}^{n}%
}\right)  $ defined on it such that:

$i)$ the laws of $(X_{\Delta_{n}},X,B,P_{X_{\Delta_{n}}},P_{X})$ and $\left(
\widehat{X_{t}^{n}},\widehat{Y_{t}^{n}},\widehat{B_{t}^{n}},\widehat{\mu
_{t}^{n}},\widehat{\nu_{t}^{n}}\right)  $ coincide for every $n\in\mathbb{N}.$

$ii)$ there exists a subsequence also denoted by $\left(  \widehat{X_{t}^{n}%
},\widehat{Y_{t}^{n}},\widehat{B_{t}^{n}},\widehat{\mu_{t}^{n}},\widehat
{\nu_{t}^{n}}\right)  $ converging to $\left(  \widehat{X_{t}},\widehat{Y_{t}%
},\widehat{B_{t}},\widehat{\mu_{t}},\widehat{\nu_{t}}\right)  $ uniformly on
every finite time interval $\widehat{P}-$a.s..

It is clear that $\left(  \widehat{B}_{t}^{n},\widehat{\mathcal{F}}_{t}%
^{n}\right)  $ and $\left(  \widehat{B}_{t},\widehat{\mathcal{F}}_{t}\right)
$ are Brownian motions with respect the filtrations $\widehat{\mathcal{F}}%
_{t}^{n}=\sigma\left(  \widehat{X_{s}^{n}},\widehat{Y_{s}^{n}},\widehat
{B_{s}^{n}};s\leq t\right)  $ and $\widehat{\mathcal{F}}_{t}=\sigma\left(
\widehat{X_{s}},\widehat{Y_{s}},\widehat{B_{s}};s\leq t\right)  .$

Note that the probability measures do not depend upon the random element
$\omega,$ then $\left(  P_{X^{n}},P_{X}\right)  =\left(  \widehat{\mu_{t}^{n}%
},\widehat{\nu_{t}^{n}}\right)  $ and consequently $\left(  \widehat{\mu
_{t}^{n}},\widehat{\nu_{t}^{n}}\right)  =\left(  P_{\widehat{X_{t}^{n}}%
},P_{\widehat{Y_{t}^{n}}}\right)  $ and $\left(  \widehat{\mu_{t}}%
,\widehat{\nu_{t}}\right)  =\left(  P_{\widehat{X_{t}}},P_{\widehat{Y_{t}}%
}\right)  $

\strut According to property $i$) and the fact that $X_{\Delta_{n}}$ and
$X_{t}$ satisfy equation (\ref{MVSDE}) and using the fact that the
finite-dimensional distributions coincide, we can easily prove that $\forall
n\geq1$ , $\forall t\geq0$
\[
E\left\vert \widehat{X_{t}^{n}}-x-%
{\displaystyle\int_{0}^{t}}
\sigma\left(  \phi_{\Delta_{n}}(s),\widehat{X_{s}^{n}},P_{\widehat{X_{t}^{n}}%
}\right)  d\widehat{B_{s}^{n}}-%
{\displaystyle\int_{0}^{t}}
b\left(  \phi_{\Delta_{n}}(s),\widehat{X_{s}^{n}},P_{\widehat{X_{t}^{n}}%
}\right)  ds\right\vert ^{2}=0,
\]

which means that%
\[
\widehat{X_{t}^{n}}=x+%
{\displaystyle\int_{0}^{t}}
\sigma\left(  \phi_{\Delta_{n}}(s),\widehat{X_{s}^{n}},P_{\widehat{X_{t}^{n}}%
}\right)  d\widehat{B_{s}^{n}}+%
{\displaystyle\int_{0}^{t}}
b\left(  \phi_{\Delta_{n}}(s),\widehat{X_{s}^{n}},P_{\widehat{X_{t}^{n}}%
}\right)  ds.
\]

Using similar arguments for $\widehat{Y_{t}^{n}}$, we obtain:%
\[
\widehat{Y_{t}^{n}}=x+%
{\displaystyle\int_{0}^{t}}
\sigma\left(  s,\widehat{Y_{t}^{n}},P_{\widehat{Y_{t}^{n}}}\right)
d\widehat{B_{s}^{n}}+%
{\displaystyle\int_{0}^{t}}
b\left(  s,\widehat{Y_{t}^{n}},P_{\widehat{Y_{t}^{n}}}\right)  ds
\]

Now, by Skorokhod's limit Theorem (see \cite{Sk} or \cite{GK} Lemma 3.1 ) and
according to $ii)$ and the fact that $\phi_{\Delta}(s)\rightarrow s,$ it holds that,%

\[%
{\displaystyle\int_{0}^{t}}
\sigma\left(  \phi_{\Delta_{n}}(s),\widehat{X}_{s}^{n_{k}},P_{\widehat{X}%
_{s}^{n_{k}}}\right)  d\widehat{B_{s}^{n_{k}}}\underset{k\rightarrow\infty
}{\overset{P}{\longrightarrow}}%
{\displaystyle\int_{0}^{t}}
\sigma\left(  s,\widehat{X}_{s},P_{\widehat{X}_{s}}\right)  d\widehat{B}_{s},
\]%
\[%
{\displaystyle\int_{0}^{t}}
b\left(  \phi_{\Delta_{n}}(s),\widehat{X}_{s}^{n_{k}},P_{\widehat{X}%
_{s}^{n_{k}}}\right)  ds\underset{k\rightarrow+\infty}{\overset{P}%
{\longrightarrow}}%
{\displaystyle\int_{0}^{t}}
b\left(  s,\widehat{X}_{s},P_{\widehat{X}_{s}}\right)  ds.
\]

We conclude that\ $\widehat{X_{t}}$ and $\widehat{Y_{t}}$ satisfy the same
stochastic differential equation (\ref{MVSDE}) on the new probability space
$\left(  \widehat{\Omega},\widehat{\mathcal{F}},\widehat{P}\right)  ,$ with
the same initial condition $x$ and common Brownian motion $\widehat{B_{t}}$ .
Therefore, according to the pathwise uniqueness for (\ref{MVSDE}) it holds
that $\widehat{X}=\widehat{Y}.$

By uniform integrability, it holds that:

$\delta\leq\underset{n\in\mathbf{N}}{\lim\inf}E\left[  \underset{t\leq T}%
{\sup}\left\vert X_{\Delta_{n}}-X_{t}\right\vert ^{2}\right]  \leq\lim
\widehat{E}\left[  \underset{t\leq T}{\sup}\left\vert \widehat{X_{t}^{n_{k}}%
}-\widehat{Y_{t}^{n_{k}}}\right\vert ^{2}\right]  =\widehat{E}\left[
\underset{t\leq T}{\sup}\left\vert \widehat{X_{t}}-\widehat{Y_{t}}\right\vert
^{2}\right]  =0$

which contradicts our hypothesis (\ref{Hyp1})$.$\eop

2) Let $(W,\mathcal{B(}W\mathcal{)},P^{B},B(t))$ be the standard Wiener
process and $X_{\Delta_{n}}(x,.,w)$ be the polygonal approximation . It is
clear that the functional $F_{\Delta_{n}}:\mathbb{R}^{d}\times W_{0}%
^{d}\longrightarrow W^{d}$ defined by $F_{_{\Delta_{n}}}(x,w)=X_{\Delta_{n}%
}(x,.,w)$ is measurable. Moreover property 1) and Borel Cantelli lemma imply
that $\left(  F_{\Delta_{n}}(x,w)\right)  $ converges uniformly in $W^{d}$ a.s..

Let $F(x,w)=\lim F_{\Delta_{n}}(x,w)$, then $F(x,w)$ is measurable and that
the unique solution is written as $X(X(0),t)=F(X(0),w)$ which achieves the
proof.\eop

\begin{remark}
1) Under the same assumptions and using the same proof, we can prove
\end{remark}

$\underset{n\rightarrow0}{\lim}\underset{}{\text{ }\sup\limits_{x\in
K}E\left[  \underset{t\leq T}{\sup}\left\vert X_{\Delta^{n}}%
(x,t)-X(x,t)\right\vert ^{2}\right]  }=0$ where $K$ is any compact set in
$\mathbb{R}^{d}.$

\section{Stability properties of MVSDEs under pathwise uniqueness}

\subsection{Stability with respect to initial conditions and coefficients}

In this section we will prove that under minimal assumptions on the
coefficients and pathwise uniqueness of solutions, the unique solution is
continuous with respect the initial condition and coefficients.

We denote by $\left(  X_{t}^{x}\right)  $ the unique solution of (\ref{MVSDE})
corresponding to the initial condition $X_{0}^{x}=x.$%

\[
\left\{
\begin{array}
[c]{l}%
dX_{t}^{x}=b(t,X_{t}^{x},\mathbb{P}_{X_{t}^{x}})dt+\sigma(t,X_{t}%
^{x},\mathbb{P}_{X_{t}^{x}})dB_{t}\\
X_{0}^{x}=x.
\end{array}
\right.
\]

\begin{theorem}
Assume that $b(t,x,\mu)$ and $\sigma(t,x,\mu)$ satisfy $\mathbf{(H_{1})}$,
$\mathbf{(H_{2}).}$Then if the pathwise uniqueness holds for equation
(\ref{MVSDE}) then the mapping
\end{theorem}

\begin{center}
$\Phi:\mathbb{R}^{d}\longrightarrow L^{2}(\Omega,\mathcal{C}(\left[
0,T\right]  ,\mathbb{R}^{d}))$
\end{center}

defined by $\left(  \Phi(x)_{t}\right)  =\left(  X_{t}^{x}\right)  $ is
continuous.\textit{\ }

\bop Suppose that the conclusion of our theorem is false, then there exists a
sequence $\left(  x_{n}\right)  $ in $\mathbb{R}^{d}$converging to $x$ and
$\delta\geq0$ such that%

\begin{equation}
\underset{n\rightarrow\infty}{\lim\inf}E\left[  \underset{t\leq T}{\sup
}\left\vert X_{t}^{n}-X_{t}\right\vert ^{2}\right]  \geq\delta\label{Hyp}%
\end{equation}

where $X_{t}^{n}=X_{t}^{x_{n}}$ and $X_{t}=X_{t}^{x}$.

Using assumptions $\mathbf{(H_{1})}$, $\mathbf{(H_{2})}$ and classical
arguments of stochastic calculus, it is easy to see that%

\[
E\left[  \left\vert X^{n}\left(  t\right)  -X^{n}\left(  s\right)  \right\vert
^{4}\right]  \leq C(T)|t-s|^{2}.
\]

where $C(T)$ is a constant which does not depend on $n$. Similar to estimate
holds true also for $X$ and the Brownian motion $B$. Then by Prokhorov's
Theorem , the sequence $(X^{n},X,P_{X^{n}},P_{X},B)$ satisfy i) and ii) of
Lemma 1 (appendix), then this sequence is tight, which implies that it is
relatively compact in the topology of weak convergence of probability
measures$.$Therefore by Skorokhod selection Theorem, there exists a
probability space $\left(  \widehat{\Omega},\widehat{\mathcal{F}},\widehat
{P}\right)  $ carrying a sequence of stochastic processes $\left(
\widehat{X_{t}^{n}},\widehat{Y_{t}^{n}},\widehat{B_{t}^{n}},\widehat{\mu
_{t}^{n}},\widehat{\nu_{t}^{n}}\right)  $ defined on it such that:

$i)$ the laws of $(X^{n},X,B,P_{X^{n}},P_{X})$ and $\left(  \widehat{X_{t}%
^{n}},\widehat{Y_{t}^{n}},\widehat{B_{t}^{n}},\widehat{\mu_{t}^{n}}%
,\widehat{\nu_{t}^{n}}\right)  $ coincide for every $n\in\mathbb{N}.$

$ii)$ there exists a subsequence also denoted by $\left(  \widehat{X_{t}^{n}%
},\widehat{Y_{t}^{n}},\widehat{B_{t}^{n}},\widehat{\mu_{t}^{n}},\widehat
{\nu_{t}^{n}}\right)  $ converging to $\left(  \widehat{X_{t}},\widehat{Y_{t}%
},\widehat{B_{t}},\widehat{\mu_{t}},\widehat{\nu_{t}}\right)  $ uniformly on
every finite time interval $\widehat{P}$-a.s., where $\left(  \widehat{B}%
_{t}^{n},\widehat{\mathcal{F}}_{t}^{n}\right)  $ and $\left(  \widehat{B}%
_{t},\widehat{\mathcal{F}}_{t}\right)  $ are Brownian motions with respect the
filtrations $\widehat{\mathcal{F}}_{t}^{n}=\sigma\left(  \widehat{X_{s}^{n}%
},\widehat{Y_{s}^{n}},\widehat{B_{s}^{n}};s\leq t\right)  $ and $\widehat
{\mathcal{F}}_{t}=\sigma\left(  \widehat{X_{s}},\widehat{Y_{s}},\widehat
{B_{s}};s\leq t\right)  .$

Note that the probability measures do not depend upon the random element
$\omega,$ then $\left(  P_{X^{n}},P_{X}\right)  =\left(  \widehat{\mu_{t}^{n}%
},\widehat{\nu_{t}^{n}}\right)  $ and consequently $\left(  \widehat{\mu
_{t}^{n}},\widehat{\nu_{t}^{n}}\right)  =\left(  P_{\widehat{X_{t}^{n}}%
},P_{\widehat{Y_{t}^{n}}}\right)  $ and $\left(  \widehat{\mu_{t}}%
,\widehat{\nu_{t}}\right)  =\left(  P_{\widehat{X_{t}}},P_{\widehat{Y_{t}}%
}\right)  $

\strut According to property $i$) and the fact that $X_{t}^{n}$ and $X_{t}$
satisfy equation (\ref{MVSDE}) with initial data $x_{n}$ and $x$, and using
the fact that the finite-dimensional distributions coincide, we can easily
prove that $\forall n\geq1$ , $\forall t\geq0$
\[
E\left\vert \widehat{X_{t}^{n}}-x_{n}-%
{\displaystyle\int_{0}^{t}}
\sigma\left(  s,\widehat{X_{s}^{n}},P_{\widehat{X_{t}^{n}}}\right)
d\widehat{B_{s}^{n}}-%
{\displaystyle\int_{0}^{t}}
b\left(  s,\widehat{X_{s}^{n}},P_{\widehat{X_{t}^{n}}}\right)  ds\right\vert
^{2}=0.
\]

In other words,%
\[
\widehat{X_{t}^{n}}=x_{n}+%
{\displaystyle\int_{0}^{t}}
\sigma\left(  s,\widehat{X_{s}^{n}},P_{\widehat{X_{t}^{n}}}\right)
d\widehat{B_{s}^{n}}+%
{\displaystyle\int_{0}^{t}}
b\left(  s,\widehat{X_{s}^{n}},P_{\widehat{X_{t}^{n}}}\right)  ds
\]

Using similar arguments for $\widehat{Y_{t}^{n}}$, we obtain:%
\[
\widehat{Y_{t}^{n}}=x+%
{\displaystyle\int_{0}^{t}}
\sigma\left(  s,\widehat{Y_{t}^{n}},P_{\widehat{Y_{t}^{n}}}\right)
d\widehat{B_{s}^{n}}+%
{\displaystyle\int_{0}^{t}}
b\left(  s,\widehat{Y_{t}^{n}},P_{\widehat{Y_{t}^{n}}}\right)  ds
\]

Now, by Skorokhod's limit theorem (see \cite{Sk} or \cite{GK} Lemma 3.1) and
according to $ii)$ it holds that,%

\[%
{\displaystyle\int_{0}^{t}}
\sigma\left(  s,\widehat{X}_{s}^{n_{k}},P_{\widehat{X}_{s}^{n_{k}}}\right)
d\widehat{B_{s}^{n_{k}}}\underset{k\rightarrow\infty}{\overset{P}%
{\longrightarrow}}%
{\displaystyle\int_{0}^{t}}
\sigma\left(  s,\widehat{X}_{s},P_{\widehat{X}_{s}}\right)  d\widehat{B}_{s},
\]%
\[%
{\displaystyle\int_{0}^{t}}
b\left(  s,\widehat{X}_{s}^{n_{k}},P_{\widehat{X}_{s}^{n_{k}}}\right)
ds\underset{k\rightarrow+\infty}{\overset{P}{\longrightarrow}}%
{\displaystyle\int_{0}^{t}}
b\left(  s,\widehat{X}_{s},P_{\widehat{X}_{s}}\right)  ds.
\]

We conclude that\ $\widehat{X_{t}}$ and $\widehat{Y_{t}}$ satisfy the same
stochastic differential equation (\ref{MVSDE}) on the new probability space
$\left(  \widehat{\Omega},\widehat{\mathcal{F}},\widehat{P}\right)  ,$ with
the same initial condition $x$ and common Brownian motion $\widehat{B_{t}}$
\[
\widehat{X_{t}}=x+%
{\displaystyle\int_{0}^{t}}
\sigma\left(  s,\widehat{X_{s}},P_{\widehat{X_{t}}}\right)  d\widehat{B_{s}}+%
{\displaystyle\int_{0}^{t}}
b\left(  s,\widehat{X_{s}},P_{\widehat{X_{t}}}\right)  ds
\]
and%

\[
\widehat{Y_{t}}=x+%
{\displaystyle\int_{0}^{t}}
\sigma\left(  s,\widehat{Y_{t}},P_{\widehat{Y_{t}}}\right)  d\widehat{B_{s}}+%
{\displaystyle\int_{0}^{t}}
b\left(  s,\widehat{Y_{t}},P_{\widehat{Y_{t}}}\right)  ds.
\]
According to the pathwise uniqueness for (\ref{MVSDE}) it holds that
$\widehat{X}=\widehat{Y}.$

By uniform integrability, it holds that:

$\delta\leq\underset{n\in\mathbf{N}}{\lim\inf}E\left[  \underset{t\leq T}%
{\sup}\left\vert X_{t}^{n}-X_{t}\right\vert ^{2}\right]  \leq\lim\widehat
{E}\left[  \underset{t\leq T}{\sup}\left\vert \widehat{X_{t}^{n_{k}}}%
-\widehat{Y_{t}^{n_{k}}}\right\vert ^{2}\right]  =\widehat{E}\left[
\underset{t\leq T}{\sup}\left\vert \widehat{X_{t}}-\widehat{Y_{t}}\right\vert
^{2}\right]  =0$

which contradicts our hypothesis (\ref{Hyp})$.$\eop

Using the same techniques we can prove the continuity of the solution of MVSDE
with respect to a parameter. In particular the solution is continunous with
respect to the coefficients. Let us consider a sequence of functions and
consider the MVSDE.

\begin{center}%
\begin{equation}
\left\{
\begin{array}
[c]{l}%
dX_{t}^{n}=\sigma_{n}\left(  t,X_{t}^{n},P_{X_{t}^{n}}\right)  dB_{t}%
+b_{n}\left(  t,X_{t}^{n},P_{X_{t}^{n}}\right)  dt\\
X^{n}(0)=x_{n}\text{.}%
\end{array}
\right.  \label{MVSDEn}%
\end{equation}
$\qquad$
\end{center}

\begin{theorem}
\textit{\ Suppose that }$\sigma_{n}\left(  t,x,\mu\right)  $\textit{\ and
}$b_{n}\left(  t,x,\mu\right)  $\textit{\ are continuous functions . Further
suppose that for each }$T>0$\textit{, and each compact set }$K$\textit{\ there
exists }$L>0$\textit{\ such that }
\end{theorem}

$\mathit{\ i)}\underset{t\leq T}{\sup}\left(  \left\vert \sigma_{n}\left(
t,x,\mu\right)  \right\vert +\left\vert b_{n}\left(  t,x,\mu\right)
\right\vert \right)  \mathit{\leq L\ }\left(  1+\left\vert x\right\vert
\right)  $ uniformly in $n$,

\textit{ii)}$\underset{n\rightarrow+\infty}{\lim}$\textit{\ }$\sup
\limits_{x\in K}\underset{t\leq T}{\sup}\left(  \left\vert \sigma_{n}\left(
t,x,\mu\right)  -\sigma\left(  t,x,\mu\right)  \right\vert +\left\vert
b_{n}\left(  t,x,\mu\right)  -b\left(  t,x,\mu\right)  \right\vert \right)
=0$,

\textit{iii) }$\lim_{n\rightarrow+\infty}x_{n}=x.$

\textit{If the pathwise uniqueness holds for equation }$(\ref{MVSDE}%
)$\textit{, then:}

\begin{center}
\textit{\ }$\underset{x\in K}{\sup}E\left[  \underset{t\leq T}{\sup}\left\vert
X_{t}^{n}-X_{t}\right\vert ^{2}\right]  =0$\textit{\ , for every }$T\geq
0$\textit{\ .}
\end{center}

\bop Similar to the proof of Theorem 4.1.\eop

\begin{remark}
\textbf{ }It is clear that if we suppose that the coefficients $b$ and
$\sigma$ are globally Lipschitz or continuous and satisfy any assumption
ensuring pathwise uniqueness, then the conclusion of the last theorem remains
true. In particular, if the coefficients satisfy Osgood condition or are
monotone ( see $\cite{BMMA, DEG})$ then the solution depends continuously on
the initial data and the coefficients $b$ and $\sigma$.
\end{remark}

\subsection{Convergence of the Picard scheme}

Assume \textbf{H}$_{1}$ and \textbf{H}$_{2}$ and consider the McKean-Vlasov
equation (\ref{MVSDE}). The sequence of successive approximations associated
to (\ref{MVSDE}) is defined by

\begin{center}%
\begin{equation}
\left\{
\begin{array}
[c]{l}%
dX_{t}^{n+1}=\sigma\left(  t,X_{t}^{n},\mathbb{P}_{X_{t}^{n}}\right)
dB_{t}+b\left(  t,X_{t}^{n},\mathbb{P}_{X_{t}^{n}}\right)  dt\\
X^{0}=x\text{.}%
\end{array}
\right.  \label{Pic(n)}%
\end{equation}
$\qquad\qquad$
\end{center}

If we assume that the coefficients are Lipschitz continuous, then the sequence
$\left(  X^{n}\right)  $ converges in quadratic mean to the unique strong
solution. In particular this method gives an effective way for the
construction of the unique strong solution $X$ of equation (\ref{MVSDE}) (see
\cite{BMMA}). Now, assume that the Lipshitz condition is not fullfilled and
assume instead only that equation (\ref{MVSDE}) admits a unique strong
solution. Does the sequence $\left(  X^{n}\right)  $ converges to $X?$ The
answer is negative even in the deterministic case. The following theorem gives
an additional necessary and sufficient condition which ensures the convergence
of successive approximations.

\begin{theorem}
\textit{\ }Assume \textbf{H}$_{1}$ and \textbf{H}$_{2}$ and that the pathwise
uniqueness\textit{ holds for (\ref{MVSDE}). Then the sequence }$\left(
X^{n}\right)  $\textit{\ converges in quadratic mean to the unique solution
}$X$ \textit{of (\ref{MVSDE})\ if and only if }$(X^{n+1}-X^{n})$%
\textit{\ converges to }$0$\textit{.}
\end{theorem}

\bop Suppose that \textbf{\ }$\left(  \mathit{X}^{n+1}-\mathit{X}^{n}\right)
$ converges to $0$ and there is some $\delta>0$ such that

\begin{center}
$\underset{n}{\inf}E\left[  \underset{t\leq T}{\sup}\left\vert X_{t}^{n}%
-X_{t}\right\vert ^{2}\right]  \geq\delta$
\end{center}

\textit{By standard techniques it is easy to see that }$\left(  X_{t}%
^{n}\right)  $ satisfies:

\textit{1) For every }$p>1,$ $\underset{n}{\sup}E\left[  \underset{t\leq
T}{\sup}\left|  X_{t}^{n}\right|  ^{2p}\right]  <+\infty.$

\textit{2) For every }$T>0$\textit{\ and }$p>1$\textit{, there exists a
constant }$C$\textit{\ independant of }$n$\textit{\ such that for every }$s<t
$\textit{\ in }$\left[  0,T\right]  $\textit{, }$E\left[  \left\vert X_{t}%
^{n}-X_{s}^{n}\right\vert ^{2p}\right]  \leq C\,\left\vert t-s\right\vert
^{p}.$

Therefore the family of stochastic process $\left(  X_{t}^{n},X_{t}%
,X_{t}^{n+1},P_{X_{t}},P_{X_{t}^{n}},P_{X_{t}^{n+1}},B_{t}\right)  $ satisfies
Kolmogorov's criteria, then by Prokhorov's Theorem this sequence of stochastic
processes is tight. Therefore by Skorokhod selection theorem, there exist a
probability space $\left(  \widehat{\Omega},\widehat{\mathcal{F}},\widehat
{P}\right)  $ carrying a sequence of stochastic processes $\left(
\widehat{X_{t}^{n}},\widehat{Y_{t}^{n}},\widehat{Z_{t}^{n}},P_{\widehat
{X_{t}^{n}}},P_{\widehat{Y_{t}^{n}}},P_{\widehat{Z_{t}^{n}}},\widehat
{B_{t}^{n}}\right)  $ with the following properties:

$\alpha)$ the laws of $\left(  X_{t}^{n},X_{t},X_{t}^{n+1},P_{X_{t}}%
,P_{X_{t}^{n}},P_{X_{t}^{n+1}},B_{t}\right)  $ and $\left(  \widehat{X_{t}%
^{n}},\widehat{Y_{t}^{n}},\widehat{Z_{t}^{n}},P_{\widehat{X_{t}^{n}}%
},P_{\widehat{Y_{t}^{n}}},P_{\widehat{Z_{t}^{n}}},\widehat{B_{t}^{n}}\right)
$ coinside for every $n\in N.$

$\beta)$ there existe a subsequence $\left(  \widehat{X_{t}^{n_{k}}}%
,\widehat{Y_{t}^{n_{k}}},\widehat{Z_{t}^{n_{k}}},P_{\widehat{X_{t}^{n_{k}}}%
},P_{\widehat{Y_{t}^{n_{k}}}},P_{\widehat{Z_{t}^{n_{k}}}},\widehat
{B_{t}^{n_{k}}}\right)  $ which converges to

$\left(  \widehat{X_{t}},\widehat{Y_{t}},\widehat{Z_{t}},P_{\widehat{X_{t}}%
},P_{\widehat{Y_{t}}},P_{\widehat{Z_{t}}},\widehat{B_{t}}\right)  $ uniformly
on every finite time interval $\widehat{P}-$a.s.

Since\textbf{\ }$\left(  \mathit{X}^{n+1}-\mathit{X}^{n}\right)  $ converges
to 0 then $\widehat{X_{t}}=\widehat{Z_{t}},$ $\widehat{P}$ a.s.

Using similar techniques as in Theorem 3.1 , $\widehat{X_{t}}$ and
$\widehat{Y_{t}}$ satisfy the same MVSDE (\ref{MVSDE}) on the new probability
space $\left(  \widehat{\Omega},\widehat{\mathcal{F}},\widehat{P}\right)  $
with the same initial condition $x$ and brownian motion $\widehat{B_{t}}.$
Then by pathwise uniqueness we have $\widehat{X}=\widehat{Y}$, $\widehat{P}$ a.s.

By uniform integrability, it holds that:
\begin{align*}
\delta &  \leq\underset{k}{\lim\inf}E\left[  \underset{0\leq t\leq T}{\sup
}\left\vert X_{t}^{n_{k}}-X_{t}\right\vert ^{2}\right]  =\text{ }\underset
{k}{\lim\inf}\widehat{E}\left[  \underset{0\leq t\leq T}{\sup}\left\vert
\widehat{X}_{t}^{n_{k}}-\widehat{Y}_{t}^{n_{k}}\right\vert ^{2}\right] \\
&  =\widehat{E}\left[  \underset{0\leq t\leq T}{\sup}\left\vert \widehat
{X}_{t}-\widehat{Y}_{t}\right\vert ^{2}\right]  =0
\end{align*}
which is a contradiction. \eop

\begin{remark}
It is well known that the convergence of the Picard scheme is equivalent to
the convergence of the series $%
{\displaystyle\sum}
(X^{n+1}-X^{n}).$ In fact Theorem 4.4 provides a simpler condition for the
convergence of $\ $this series.
\end{remark}

Let us consider MVSDEs of the type (\ref{MVSDE}) whose coefficients are real
valued and satisfy:

$\mathbf{(A}_{\mathbf{1}}\mathbf{)}$ There exist $C>0,$ such that for every
$x\in\mathbb{R}$ and $\left(  \mu,\nu\right)  \in\mathcal{P}_{1}%
(\mathbb{R})\times\mathcal{P}_{1}(\mathbb{R}):$

\begin{center}
$|b(t,x,\mu)-b(t,x,\nu)|\leq CW_{2}(\mu,\nu)$
\end{center}

$\mathbf{(A}_{\mathbf{2}}\mathbf{)}$ There exist $K>0,$ such that for every
$x,y\in\mathbb{R}$ $,$ $|\sigma(t,x)-\sigma(t,y)|\leq K|x-y|).$

$\mathbf{(A}_{\mathbf{3}}\mathbf{)}$ There exists a strictly increasing
function $\kappa(u)$ on $[0,+\infty)$ such that $\kappa(0)=0$ and $\kappa$ is
concave satisfying $%
{\displaystyle\int\limits_{0^{+}}}
\kappa^{-1}(u)du=+\infty,$ such that for every $(x,y)\in\mathbb{R}^{d}%
\times\mathbb{R}^{d}$ and $\mu\in\mathcal{P}_{1}(\mathbb{R}),$ $|b(t,x,\mu
)-b(t,y,\mu)|\leq\kappa(|x-y|)$

\begin{corollary}
Assume \textbf{A}$_{\mathbf{1}}-$\textbf{A}$_{3}.$ Then equation (\ref{MVSDE})
has a unique strong solution such that the Picard scheme of successive
approximations converges.
\end{corollary}

\bop

The existence and uniqueness of a strong solution was proved in $\cite{BMMA}$.
According to Theorem 4.4, to prove the convergence of the successive
approximations it is sufficient to show that $(X^{n+1}-X^{n})$ converges to
$0$.

Let $u_{n}(t)=E\left(  \underset{0\leq s\leq t}{\sup}\left\vert X_{s}%
^{n+1}-X_{s}^{n}\right\vert \right)  $ by following the same steps as in
\cite{BMMA} Theorem 3.2 and using the concavity of the function $\kappa$, it
holds that: $u_{n+1}(t)\leq M%
{\displaystyle\int_{0}^{t}}
\kappa\left(  u_{n}(s)\right)  ds,$ $t\in\left[  0,T\right]  .$

Let $\left(  \varphi_{n}\right)  $ be a sequence of functions defined by:

\begin{center}
$\varphi_{0}(t)=y(t)$

$\varphi_{n+1}(t)=M%
{\displaystyle\int_{0}^{t}}
\kappa\left(  \varphi_{n}(s)\right)  ds$
\end{center}

By using a lemma in \cite{Di} page 114-124, it is possible to choose the
function $y(t)$ such that:

\begin{center}
$y(t)\geq u_{0}(t)$

$y(t)\geq M%
{\displaystyle\int_{0}^{t}}
\kappa\left(  y(s)\right)  ds$
\end{center}

By induction, it is easy to see that $u_{n}(s)\leq\varphi_{n}(s)$ and $\left(
\varphi_{n}\right)  $ is decreasing. Denote $\varphi=\lim\varphi_{n}.\ $Note
that the convergence is uniform and $\varphi$ is continuous and satisfies:

\begin{center}
$\varphi(t)\geq M%
{\displaystyle\int_{0}^{t}}
\kappa\left(  \varphi(s)\right)  ds.$

\end{center}

Assumption A3) implies that $\varphi=0.$ Therefore $\lim u_{n}(t)=0$, which
achieves the proof$.$\eop

\subsection{Stability of MVSDEs driven by semi-martingales}

In this section, we consider the McKean-Vlasov equations driven by continuous semi-martingales.

Let $b:\left[  0,1\right]  \times\mathbf{R}^{d}\times\mathcal{P}%
_{2}(\mathbf{R}^{d})\longrightarrow\mathbf{R}^{d}$ and $\sigma:\left[
0,1\right]  \times\mathbf{R}^{d}\times\mathcal{P}_{2}(\mathbf{R}%
^{d})\longrightarrow\mathbf{R}^{d\times d}$ be bounded continuous functions.

We consider the stochastic differential equation:%

\begin{equation}
\left\{
\begin{array}
[c]{l}%
dX_{t}=\sigma\left(  t,X_{t},\mathbb{P}_{X_{t}}\right)  dM_{t}+b\left(
t,X_{t},\mathbb{P}_{X_{t}}\right)  dA_{t}\\
X_{0}=x
\end{array}
\right.  \label{Mart}%
\end{equation}
\newline where $A_{t}$ is an adapted continuous process of bounded variation
and $M_{t}$ is a continuous local martingale.

\begin{definition}
\textit{\ Pathwise uniqueness property holds for equation \ref{Mart}, if
whenever }$\left(  X,M,A,\left(  \Omega,\mathcal{F}\text{,}P\right)
,\mathcal{F}_{t}\right)  $\textit{\ and }$\left(  X^{\prime},M^{\prime
},A^{\prime},\left(  \Omega,\mathcal{F}\text{,}P\right)  ,\mathcal{F}%
_{t}^{\prime}\right)  $\textit{\ are two weak solutions such that }$\left(
M,A\right)  =\left(  M^{\prime},A^{\prime}\right)  $\textit{\ P a.s, then
}$X=X^{\prime}$\textit{\ P- a.s.}
\end{definition}

We consider the following equations:

\begin{center}%
\begin{equation}
\left\{
\begin{array}
[c]{l}%
dX_{t}^{n}=\sigma\left(  t,X_{t}^{n}\right)  dM_{t}^{n}+b\left(  t,X_{t}%
^{n}\right)  dA_{t}^{n}\\
X_{0}^{n}=x\text{.}%
\end{array}
\right.  \label{Mart(n)}%
\end{equation}
$\qquad\qquad$
\end{center}

where $\left(  M^{n}\right)  $ is a sequence of continuous $\left(
\mathcal{F}_{t},P\right)  -$local martingales and $\left(  A^{n}\right)  $ a
sequence of $\mathcal{F}_{t}$-adapted continuous processes with bounded variation.

Let us suppose that $\left(  A,A^{n},M,M^{n}\right)  $ satisfy the following conditions:

$(\mathbf{H}_{4})$ The family $\left(  A,A^{n},M,M^{n}\right)  $ is bounded in
probability in $C\left(  \left[  0,1\right]  \right)  ^{4}.$

$(\mathbf{H}_{5})$ $M^{n}-M$ $\longrightarrow0$ in probability in $C\left(
\left[  0,1\right]  \right)  .$

$(\mathbf{H}_{6})$ The total variation $\left(  A^{n}-A\right)
\longrightarrow0$ in probability as $n\rightarrow+\infty.$

In the proof of the main result of this section, we need the following lemmas
given in \cite{Gy} which are variants of Skorokhod's limit theorem( \cite{Sk}
or \cite{GK} Lemma 3.1).

\begin{lemma}
\textit{Let }$\left\{  f_{n}\left(  t\right)  ,f\left(  t\right)  \right\}
$\textit{\ be a family of continuous processes and let }$\left\{  C_{n}\left(
t\right)  ,C\left(  t\right)  \right\}  $\textit{\ be a family of continuous
processes of bounded variation. Assume that: }
\end{lemma}

\begin{center}
\textit{\ }$\underset{n\rightarrow+\infty}{\lim}f_{n}=f$\textit{\ in
probability in }$C\left(  \left[  0,T\right]  \right)  $\textit{.}

\textit{\ }$\underset{n\rightarrow+\infty}{\lim}C_{n}=C$\textit{\ in
probability in }$C\left(  \left[  0,T\right]  \right)  $\textit{.}

$\left\{  Var\left(  C_{n}\right)  \text{ ; }n\in\mathbf{N}\right\}
$\textit{\ is bounded in probability.}
\end{center}

\textit{Then the following result holds:}

\begin{center}
$\forall\varepsilon>0$\textit{,}$\qquad\underset{n\rightarrow+\infty}{\lim
}P\left[  \underset{t\leq1}{\sup}\left\vert
{\displaystyle\int_{0}^{T}}
f_{n}dC_{n}-%
{\displaystyle\int_{0}^{T}}
fdC\right\vert >\varepsilon\right]  =0$.
\end{center}

\begin{lemma}
\textit{Consider a family of filtrations }$\left(  F_{t}^{n}\right)  ,\left(
F_{t}\right)  $\textit{\ satisfying the usual conditions. Let }$\left\{
f_{n}\left(  t\right)  ,f\left(  t\right)  :t\in\left[  0,T\right]  \right\}
$\textit{\ be a sequence of continuous adapted processes and let }$\left\{
N_{n}\left(  t\right)  ,N\left(  t\right)  :t\in\left[  0,T\right]  \right\}
$\textit{\ be a sequence of continuous local martingales with respect to
}$\left(  F_{t}^{n}\right)  ,\left(  F_{t}\right)  $\textit{\ respectively.
Suppose that}
\end{lemma}

\begin{center}
$\underset{n\rightarrow+\infty}{\lim}f_{n}=f$\textit{\ in probability in
}$C\left(  \left[  0,T\right]  \right)  .$

$\underset{n\rightarrow+\infty}{\lim}N_{n}=N$\textit{\ in probability in
}$C\left(  \left[  0,T\right]  \right)  .$
\end{center}

\textit{Then}

\begin{center}
$\forall\varepsilon>0$\textit{\ ,}$\qquad\underset{n\rightarrow+\infty}{\lim
}P\left[  \underset{t\leq1}{\sup}\left\vert
{\displaystyle\int_{0}^{T}}
f_{n}dN_{n}-%
{\displaystyle\int_{0}^{T}}
fdN\right\vert >\varepsilon\right]  =0$.
\end{center}

\begin{theorem}
\textit{\ Assume }$\mathbf{H}_{1}-\mathbf{H}_{6}$\textit{\ . If the pathwise
uniqueness holds for equation \ref{Mart} then:}
\end{theorem}

\begin{center}%
\[
\lim_{n\rightarrow\infty}E[\sup_{t\leq T}|X_{t}^{n}-X_{t}|^{2}]=0
\]
\smallskip\ 
\end{center}

\bop

Suppose that the conclusion of our theorem is false. Then there exists
$\delta>0$ such that%

\[
\underset{n}{\inf}E[\sup_{t\leq T}|X_{t}^{n}-X_{t}|^{2}]\geq\delta
\]

The assumptions made on the coefficients and the driving processes imply that
$X^{n}$ and $X$ are tight. Then, the family $\Gamma^{n}=\left(  X^{n}%
,X,\mathbb{P}_{X^{n}},\mathbb{P}_{X},A^{n},A,M^{n},M\right)  $ is tight.
Therefore by Skorokhod's selection Theorem there exist a probability space
$\left(  \Omega^{\prime},\mathcal{F}^{\prime}\text{,}P^{\prime}\right)  $
carrying a sequence $\Gamma^{\prime n}=\left(  X^{\prime n},Y^{\prime
n},\mathbb{P}_{X^{\prime n}},\mathbb{P}_{Y^{\prime n}},A^{\prime n},B^{\prime
n},M^{\prime n},N^{\prime n}\right)  $ such that:

i) For each $n,$ $\Gamma^{n}$ and $\Gamma^{\prime n}$ have the same distribution.

ii) There exists a subsequence $\left(  \Gamma^{\prime n_{k}}\right)  $ still
denoted by $\left(  \Gamma^{\prime n}\right)  $ which converges $P^{\prime}$
a.s to $\Gamma^{\prime}$ where

$\Gamma^{\prime n}=\left(  X^{\prime n},Y^{\prime n},\mathbb{P}_{X^{\prime n}%
},\mathbb{P}_{Y^{\prime n}},A^{\prime},B^{\prime},M^{\prime},N^{\prime
}\right)  .$

Let $\left(  \mathcal{F}_{t}^{\prime n};t\in\left[  0,T\right]  \right)  $ be
the filtration generated by $Z_{t}^{\prime n}$ $(t\in\left[  0,T\right]  )$
and $\left(  \mathcal{F}_{t}^{\prime};t\in\left[  0,T\right]  \right)  $ be
the filtration generated by $Z_{t}^{\prime}$ $(t\in\left[  0,T\right]  )$
then, $M_{t}^{\prime n},N^{\prime n}$ (resp. $M_{t}^{\prime},N^{\prime}$) are
$\mathcal{F}_{t}^{\prime n}$ (resp. $\mathcal{F}_{t}^{\prime}$) continuous
local martingales. The processes $X^{\prime n}$ and $Y^{\prime n}$ satisfy the
following equations :

\begin{center}%
\[
\left\{
\begin{array}
[c]{l}%
dX_{t}^{\prime n}=\sigma\left(  t,X_{t}^{\prime n},\mathbb{P}_{X_{t}^{\prime
n}}\right)  dM_{t}^{\prime n}+b\left(  t,X_{t}^{\prime n},\mathbb{P}%
_{X_{t}^{\prime n}}\right)  dA_{t}^{\prime n}\\
X_{0}^{\prime n}=x
\end{array}
\right.
\]
$\qquad\qquad$%

\[
\left\{
\begin{array}
[c]{l}%
dY_{t}^{\prime n}=\sigma\left(  t,Y_{t}^{\prime n},\mathbb{P}_{Y_{t}^{\prime
n}}\right)  dN^{\prime n}+b\left(  t,Y_{t}^{\prime n},\mathbb{P}%
_{Y_{t}^{\prime n}}\right)  dB^{\prime n}\\
Y_{0}^{\prime n}=x
\end{array}
\right.
\]
$\qquad\qquad$
\end{center}

By using Lemmas 4.6 and 4.7, we see that the limiting processes satisfy the
following equations:

\begin{center}
$\left\{
\begin{array}
[c]{l}%
dX_{t}^{\prime}=\sigma\left(  t,X_{t}^{\prime},\mathbb{P}_{X_{t}^{\prime}%
}\right)  dM_{t}^{\prime}+b\left(  t,X_{t}^{\prime},\mathbb{P}_{X_{t}^{\prime
}}\right)  dA_{t}^{\prime}\\
X_{0}^{\prime}=x
\end{array}
\right.  \qquad\qquad$

$\left\{
\begin{array}
[c]{l}%
dY_{t}^{\prime}=\sigma\left(  t,Y_{t}^{\prime},\mathbb{P}_{Y_{t}^{\prime}%
}\right)  dN_{t}^{\prime}+b\left(  t,Y_{t}^{\prime},\mathbb{P}_{Y_{t}^{\prime
}}\right)  dB_{t}^{\prime}\\
Y_{0}^{\prime}=x\text{.}%
\end{array}
\right.  \qquad\qquad$
\end{center}

By using hypothesis $\left(  \mathbf{H}_{4}\right)  $ and $\left(
\mathbf{H}_{5}\right)  $, it is easy to see that $M^{\prime}=\widetilde{M}$
and $A^{\prime}=\widetilde{A}$, $P^{\prime}$ $a.s$.

Hence by pathwise uniqueness $X^{\prime}=Y^{\prime}$. This contradicts our
assumption, therefore

\begin{center}
$\underset{n\rightarrow\infty}{\lim}E\left[  \underset{t\leq T}{\sup
}\left\vert X_{t}^{n}-X_{t}\right\vert ^{2}\right]  =0.$\eop

\end{center}

\section{Existence and uniqueness is a generic property}

We know that under globally Lipschitz coefficients equation (\ref{MVSDE}) has
a unique strong solution (see \cite{Gra, JMW, Sn}). A huge literature has been
produced to improve the conditions under which pathwise uniquenes holds.
Moreover the continuity of the coefficients is not sufficient for the
uniqueness. The objective to identify completely the set of coefficients,
under which there is a unique strong solution seems to be out of reach, even
for ordinary differential equations. In this section we are interested in
qualitative properties of the set of coefficients for which existence an
uniqueness of solutions hold. In fact we prove that "most" of the MVSDEs with
bounded uniformly continuous coefficients enjoy the property of existence and
uniqueness. The expression "most" should be understood in the sense of
topology and is similar to the measure theoritic concept of a set whose
complement is a negligible set. More precisely, we prove that in the sense of
Baire, the set of coefficients $\left(  b,\sigma\right)  $ for which existence
and uniqueness of a strong solution is a residual set in the Baire space of
all bounded uniformly continuous functions.

Prevalence properties for ordinary differential equations were first
considered by Orlicz \cite{Or} and Lasota-Yorke \cite{LaYo}$.$ These
properties have been extended to It\^{o} stochastic differential equations in
\cite{BMO1, BMO2, Heu}$.$

Let us recall some facts about Baire spaces.

\begin{definition}
A Baire space $X$ is a topological space in which the union of every countable
collection of closed sets with empty interior has empty interior.
\end{definition}

This definition is equivalent to each of the following conditions.

a) Every intersection of countably many dense open sets is dense.

b) The interior of every union of countably many closed nowhere dense sets is empty.

\begin{remark}
By the Baire category theorem, we know that a complete metric space is a Baire space.
\end{remark}

\begin{definition}
1) A subset of a topological space $X$ is called nowhere dense in $X,$ if the
interior of its closure is empty

2) A subset is of first category in the sense of Baire (or meager in X), if it
is a union of countably many nowhere dense subsets.

3) A subset is of second category or nonmeager in $X$, if it is not of first
category in $X.$
\end{definition}

\begin{remark}
1) The definition for a Baire space can then be stated as follows: a
topological space X is a Baire space if every non-empty open set is of second
category in X.

2) In the literature, a subset of second category is also called a residual
subset. \ 
\end{remark}

\begin{definition}
A property $P$ is generic in the Baire space $\mathcal{X}$ if $P$ holds is
satisfied for each element in $\mathcal{X}-\mathcal{N}$, where $\mathcal{N}$
is a set of first category in the Baire space $\mathcal{X}$.
\end{definition}

Let us introduce the appropriate Baire space.

Let $\mathcal{C}_{1}$be the set of bounded uniformly continuous functions
$b:\mathbb{R}_{+}\times\mathbb{R}^{d}\times\mathcal{P}_{2}\left(
\mathbb{R}^{d}\right)  \longrightarrow\mathbb{R}^{d}$. Define the metric
$\rho_{1}$ on $\ \mathcal{C}_{1}$ as follows:

\begin{center}
$\rho_{1}\left(  b_{1},b_{2}\right)  =\underset{(t,x,\mu)\in\mathbb{R}%
_{\mathbb{+}}\times\mathbb{R}^{d}\times\mathcal{P}_{2}\left(  \mathbb{R}%
^{d}\right)  }{\sup}\left\vert b_{1}(t,x,\mu)-b_{2}(t,x,\mu)\right\vert
\newline$
\end{center}

Note that the metric $\rho_{1}$ is compatible with the topology of uniform
convergence on $\mathbb{R}_{\mathbb{+}}\times\mathbb{R}^{d}\times
\mathcal{P}_{2}\left(  \mathbb{R}^{d}\right)  .$

Let $\mathcal{C}_{2}$ be the set of bounded unifomly continuous functions
$\sigma:\mathbb{R}_{+}\times\mathbb{R}^{d}\times\mathcal{P}_{2}\left(
\mathbb{R}^{d}\right)  \longrightarrow\mathbb{R}^{d}\otimes\mathbb{R}^{d}$
endowed with the corresponding metric $\rho_{2}:$

\begin{center}
$\rho_{2}\left(  \sigma_{1},\sigma_{2}\right)  =\underset{(t,x,\mu
)\in\mathbb{R}_{\mathbb{+}}\times\mathbb{R}^{d}\times\mathcal{P}_{2}\left(
\mathbb{R}^{d}\right)  }{\sup}\left\vert \sigma_{1}(t,x,\mu)-\sigma
_{2}(t,x,\mu)\right\vert $
\end{center}

It is clear that since $\mathbb{R}_{\mathbb{+}}\times\mathbb{R}^{d}%
\times\mathcal{P}_{2}\left(  \mathbb{R}^{d}\right)  $ is a complete metric
space, then $\mathfrak{R}=\mathcal{C}_{1}\times\mathcal{C}_{2}$ endowed with
the metric $\lambda$ is a complete metric space also, where $\lambda(\left(
b_{1},\sigma_{1}\right)  ,\left(  b_{2},\sigma_{2}\right)  )=\rho\left(
b_{1},b_{2}\right)  +\rho\left(  \sigma_{1},\sigma_{2}\right)  .$

\begin{remark}
\textbf{ } Note that for ordinary or It\^{o} stochastic differential
equations, the suitable Baire space is the space of bounded continous
functions. The space of continuous functions contains a dense subset formed of
all locally Lipschitz functions for which there is uniqueness of solutions for
It\^{o} SDEs. This property is no more valid for MVSDEs as the uniqueness of
solutions may fail for locally Lipschitz coefficients (see \cite{Sheu} ).
Instead of bounded continuous functions we consider bounded uniformly
continuous functions.\ These functions are approximated by globally Lipschitz
functions for which we have existence and uniquness. The fact that the
coefficients depend on the marginal distribution of the unknown process is not
suitable for localization techniques$.$
\end{remark}

For $\left(  b,\sigma\right)  $ in $\mathfrak{R,}$ let $E(x,b,\sigma)$ stands
for MVSDE (\ref{MVSDE}) corresponding to coefficients $b,\sigma$ and initial
data $x.$

\begin{center}
$\mathbf{M}^{2}=\left\{  \xi:\Omega\times\mathbb{R}_{+}\longrightarrow
\mathbb{R}^{d}\text{, }F_{t}^{B}-\text{adapted, continuous with }E\left[
\underset{t\leq T}{\sup}\left\vert \xi_{t}\right\vert ^{2}\right]
<+\infty\right\}  $
\end{center}

Define a metric on $\mathbf{M}^{2}$ by:

\begin{center}
$d\left(  \xi_{1},\xi_{2}\right)  =\left(  E\underset{0\leq t\leq T}{\sup
}\left\vert \xi_{t}^{1}-\xi_{t}^{2}\right\vert ^{2}\right)  ^{\frac{1}{2}}$
\end{center}

By using Borel-Cantelli lemma, it is easy to see that $\left(  \mathbf{M}%
^{2},d\right)  $ is a complete metric space.

It is clear that a strong solution ($\xi_{t})$ of equations (\ref{MVSDE}) is
an element of the metrix space (\textbf{M}$^{2},d).$

Let $\mathcal{L}$ be the subset of $\mathfrak{R}$ consisting of functions
$h(t,x,\mu)$ which are Lipschitz in their arguments, that is:

\begin{center}
$|b(t,x,\mu)-b(t,x^{\prime},\mu^{\prime})|\leq C\left(  \left\vert
x-x^{\prime}\right\vert +W_{2}(\mu,\mu^{\prime}\right)  ,$

$|\sigma(t,x,\mu)-\sigma(t,x^{\prime},\mu^{\prime})|\leq C\left(  \left\vert
x-x^{\prime}\right\vert +W_{2}(\mu,\mu^{\prime}\right)  $
\end{center}

\begin{proposition}
Every bounded uniformly continuous function in a metric space is a uniform
limit of a sequence of globally Lipschitz functions.
\end{proposition}

\bop See \cite{Hei} Theorem 6.8\eop

The last proposition states that the subset $\mathcal{L}$ of globally
Lipschitz functions\ is dense in the Baire space $\mathcal{R}$.

\subsection{The oscillation function}

Let us define the oscillation function, which was first introduced by
Lasota-Yorke \cite{LaYo} in the case of ordinary differential equations and
partial differential equations and then used by \cite{BMO2, Heu} for It\^{o} SDEs.

Let $x\in\mathbb{R}^{d}$ and $(b,\sigma)\in\mathfrak{R},$ let $\xi
(x,b,\sigma)$\textit{\ }the solution\textit{ }of equation $E(x,b,\sigma).$

Define the oscillation function as follows

\begin{center}
$D_{1}(x,b,\sigma):\mathbb{R}^{d}\times\mathcal{R}\longrightarrow
\mathbb{R}_{+}$
\end{center}

$D_{1}(x,b,\sigma)=\lim_{\delta\rightarrow0}\sup\left\{  d(\xi(x,b_{1}%
,\sigma_{1}),\xi(x,b_{2},\sigma_{2});\text{ }\left(  b_{i},\sigma_{i}\right)
\in\mathcal{L}\text{ and }\lambda(\left(  b,\sigma\right)  ,(b_{i},\sigma
_{i}))<\delta,\text{ }i=1,2\right\}  $

\begin{proposition}
Let $x\in\mathbb{R}^{d}$ and $\left(  b,\sigma\right)  $ are Lipschitz
coefficients, that is $\left(  b,\sigma\right)  \in\mathcal{L}$, then
$D_{1}(x,b,\sigma)=0.$
\end{proposition}

\bop

We know that if $\left(  b,\sigma\right)  \in\mathcal{L}$ then equation
(\ref{MVSDE}) has a unique strong solution.

For each $i=1,2$, let $\left(  X_{t}^{i}\right)  $ be a solution of
(\ref{MVSDE}) coerresponding to $(b_{i},\sigma_{i})$, then%

\begin{align*}
|X_{t}^{1}-X_{t}^{2}|^{2}  &  \leq3(\int_{0}^{t}|b_{1}(s,X_{s}^{1}%
,\mathbb{P}_{X_{s}^{1}})-b_{1}(s,X_{s}^{2},\mathbb{P}_{X_{s}^{2}})|ds)^{2}\\
&  +3(\int_{0}^{t}|b_{1}(s,X_{s}^{2},\mathbb{P}_{X_{s}^{2}})-b_{2}(s,X_{s}%
^{2},\mathbb{P}_{X_{s}^{2}}))|ds)^{2}\\
&  +3\left\vert \int_{0}^{t}\left(  \sigma_{1}(s,X_{s}^{1},\mathbb{P}%
_{X_{s}^{1}})-\sigma_{1}(s,X_{s}^{2},\mathbb{P}_{X_{s}^{2}})\right)
dB_{s}\right\vert ^{2}\\
&  +3\left\vert \int_{0}^{t}\left(  \sigma_{1}(s,X_{s}^{2},\mathbb{P}%
_{X_{s}^{2}})-\sigma_{2}(s,X_{s}^{2},\mathbb{P}_{X_{s}^{2}})\right)
dB_{s}\right\vert ^{2}.
\end{align*}

By using the Lipschitz continuity and Burkholder Davis Gundy inequality, it
holds that%

\begin{align*}
E\left[  \sup_{t\leq T}|X_{t}^{1}-X_{t}^{2}|^{2}\right]   &  \leq
3(T+C_{2})L^{2}\int_{0}^{t}\left(  E\left[  \sup_{s\leq t}|X_{s}^{1}-X_{s}%
^{2}|^{2}\right]  +W_{2}(\mathbb{P}_{X_{s}^{1}},\mathbb{P}_{X_{s}^{2}}%
)^{2}\right)  ds\\
&  +6(T+C_{2})E[\int_{0}^{t}|b_{1}(s,X_{s}^{2},\mathbb{P}_{X_{s}^{2}%
})-b(s,X_{s}^{2},\mathbb{P}_{X_{s}^{2}}|^{2}ds]\\
&  +6(T+C_{2})E[\int_{0}^{t}\int_{0}^{t}|\sigma_{1}(s,X_{s}^{2},\mathbb{P}%
_{X_{s}^{2}})-\sigma(s,X_{s}^{2},\mathbb{P}_{X_{s}^{2}}|^{2}ds]\\
&  \leq6(T+C_{2})\int_{0}^{t}E\left[  \sup_{s\leq t}|X_{s}^{1}-X_{s}^{2}%
|^{2}\right]  ds+K,
\end{align*}

such that
\[
K=3(T+C_{2})E[\int_{0}^{T}\left\vert b_{1}-b\right\vert ^{2}(s,X_{s}%
^{2},\mathbb{P}_{X_{s}^{2}})+\left\vert \sigma_{1}-\sigma\right\vert
^{2}(s,X_{s}^{2},\mathbb{P}_{X_{s}^{2}})ds].
\]

An application of Gronwall lemma allows us to get%

\[
E\left[  \sup_{t\leq T}|X_{t}^{1}-X_{t}^{2}|^{2}\right]  \leq C\delta^{2}.
\]

where C is some constantwhcih implies that $D_{1}(x,b,\sigma)=0$.\eop

\begin{proposition}
The oscillation function $D$ is upper semicontinuous function at each point of
the set $\mathbb{R}^{d}\times\mathcal{L}.$
\end{proposition}

\bop

Let $\left(  x_{n},b_{n},\sigma_{n}\right)  $ be a sequence in $\mathbb{R}%
^{d}\times\mathcal{R}$ converging to a limit $\left(  x,b,\sigma\right)
\in\mathbb{R}^{d}\times\mathcal{L}$. $D$ is upper semicontinuous if
$\lim_{n\rightarrow+\infty}D_{1}\left(  x_{n},b_{n},\sigma_{n}\right)  =0.$
Suppose that the last statement is false. Then according to the definition of
the function $D,$ there exists $\varepsilon>0$ and a subsequence still denoted
by $\left\{  n\right\}  $ (to avoid heavy notations) and functions $\left(
b_{n}^{i},\sigma_{n}^{i}\right)  $ in $\mathcal{L}$ such that :

\begin{center}
(i) $\lambda(\left(  b_{n},\sigma_{n}\right)  ,\left(  b_{n}^{i},\sigma
_{n}^{i}\right)  )<1/2^{n}$

(ii) $d(\xi(x_{n},b_{n}^{1},\sigma_{n}^{1}),\xi(x_{n},b_{n}^{2},\sigma_{n}%
^{2}))>\varepsilon/2$
\end{center}

Thus according to Theorem 4.1 on the continuous dependence with respect to
initial condition and coefficients, and property (i) it holds that:

\begin{center}
$\lim_{n\rightarrow+\infty}d(\xi(x_{n},b_{n}^{1},\sigma_{n}^{1}),\xi
(x_{n},b_{n}^{2},\sigma_{n}^{2}))=0.$
\end{center}

But this contradicts the property (ii), then $D$ is upper
semicontinuous.\eop

\begin{proposition}
Let $\left(  x,b,\sigma\right)  $ be in $\mathbb{R}^{d}\times\mathcal{R}$ such
that $D_{1}\left(  x,b,\sigma\right)  =0,$ then there exists at least one
strong solution to MVSDE (\ref{MVSDE})$.$
\end{proposition}

\bop Similar to \cite{BMO2} Prop. 1.4 or Proposition 5 in \cite{Heu}$.$\eop

\subsection{Existence and uniqueness of solutions is a generic propery}

The main result of this section is the following.

\begin{theorem}
\textit{The subset }$\mathcal{U}$\textit{\ \ consisting of those }$\left(
\sigma,b\right)  $\textit{\ for which existence and uniqueness of a strong
solution holds for }equation (\ref{MVSDE})\textit{\ contains a set of second
category in the Baire space }$\mathfrak{R}$\textit{.}
\end{theorem}

\bop\textbf{ }It is clear from Proposition 5.10, that if for some $\left(
\sigma,b\right)  $ in $\mathfrak{R}$, equation (\ref{MVSDE}) has at least one
strong solution then $D_{1}\left(  x,b,\sigma\right)  =0$. Then the set of
couples $\left(  \sigma,b\right)  $ in $\mathfrak{R,}$ for which existence of
strong solution holds, contains the set

\begin{center}
$\mathcal{A}=\left\{  \left(  \sigma,b\right)  \in\mathfrak{R;}D_{1}\left(
x,b,\sigma\right)  =0\right\}  .$
\end{center}

If we denote

\begin{center}
$\mathcal{A}_{n}=\left\{  \left(  \sigma,b\right)  \in\mathfrak{R;}%
D_{1}\left(  x,b,\sigma\right)  <1/n\right\}  $
\end{center}

then $\mathcal{A}=%
{\displaystyle\bigcap\limits_{n=1}^{+\infty}}
\mathcal{A}_{n}.$

Let $\left(  b,\sigma\right)  \in\mathcal{L}$, then according to Proposition
5.8 we have $D_{1}\left(  x,b,\sigma\right)  =0.$ Therefore $\mathcal{L}%
\subset\mathcal{A}_{n}$ and then by Proposition 5.7, $\mathcal{A}_{n}$
contains a dense open subset of the Baire space $\left(  \mathfrak{R}%
,\lambda\right)  .$ Therefore $\mathcal{A}$ contains an intersection of open
dense subsets, then $\mathcal{A}$ is a residual subset in $\mathfrak{R}$.

If $\left(  b,\sigma\right)  \in\mathcal{A}$, equation MVSDE (\ref{MVSDE})
enjoys the property of existence of a solution. To obtain the property of
uniqueness let us introduce the finction $D_{2}:\mathcal{A}\longrightarrow
\left[  0,+\infty\right[  $ defined by

\begin{center}
$D_{2}(\left(  b,\sigma\right)  )=\sup\left\{  d(\xi_{1},\xi_{2});\xi_{i}%
\in\mathbf{S}^{2}\text{ and }\xi_{i}\text{ is a strong solution of
}E(x,b\text{,}\sigma)\right\}  $
\end{center}

and

\begin{center}
$\mathcal{B}_{n}=\left\{  \left(  \sigma,b\right)  \in\mathcal{A}%
\mathfrak{;}D_{2}\left(  x,b,\sigma\right)  <1/n\right\}  $
\end{center}

Let $\mathcal{B}=%
{\displaystyle\bigcap\limits_{n=1}^{+\infty}}
\mathcal{B}_{n}$

Let us note that if $\left(  b,\sigma\right)  $ are Lipschitz functions then
equation (\ref{MVSDE}) admits a unique strong solution. Therefore if $\left(
b,\sigma\right)  \in\mathcal{L}$, then $D_{2}\left(  x,b,\sigma\right)  =0.$
This implies in particular that $\mathcal{B}_{n}$ contains the intersection of
$\mathcal{A}$ and a dense open subset in $\mathcal{R}$, namely $\mathcal{L}$.
Therefore $\mathcal{B}$ contains an intersection of open dense subsets in the
Baire space $\left(  \mathfrak{R},\lambda\right)  .$ This means that
$\mathcal{B}$ is a residual subset in $\left(  \mathfrak{R},\lambda\right)
.$\eop

\begin{remark}
By using similar techniques it is not difficult to prove that the set of
coefficients $\left(  b,\sigma\right)  $ for which the Euler polygonal scheme
and the Picard scheme for MVSDE converge, is a residual set in the Baire space
of all bounded uniformly continuous functions $\left(  \mathfrak{R}%
,\lambda\right)  .$
\end{remark}

\end{document}